\documentclass[a4paper,leqno]{article}
\usepackage{amsmath}
\usepackage{amsfonts}
\usepackage{amssymb}
\usepackage{amsthm}

\theoremstyle{plain}
\newtheorem{theorem}{Theorem}

\newtheorem{corollary}[theorem]{Corollary}
\newtheorem*{Theorem}{Theorem}
\newtheorem*{Lemma}{Lemma}

\newtheorem*{Problem}{Problem}

\theoremstyle{definition}

\newtheorem*{Definition}{Definition}

\theoremstyle{remark}

\title{NEW RECURRENT INEQUALITY ON A CLASS \\ OF VERTEX FOLKMAN NUMBERS}

\author{Nikolay Rangelov Kolev \and Nedyalko Dimov Nenov}

\date{}

\begin{document}

\maketitle

{
  \parindent0pt \footnotesize \leftskip20pt \rightskip20pt
  \baselineskip10pt

  Let $G$ be a graph and $V(G)$ be the vertex set of $G.$ Let $a_1$
  ,\dots , $a_r$ be positive integers, $m=\sum_{i=1}^{r} (a_{i}-1)+1$
  and $p= \max \{a_1, \ldots ,a_r\}$.  The symbol $G\rightarrow \{a_1
  , \ldots , a_r\}$ denotes that in every $r$-coloring of $V(G)$ there
  exists a monochromatic $a_i$-clique of color $i$ for some $i=1 ,
  \ldots , r$. The vertex Folkman numbers $ F(a_1 , \dots , a_r ;
  m-1)=\min \{| V(G) | : G\rightarrow (a_1\dots a_r)$ and $K_{m-1}
  \not \subseteq G \}$ are considered. In this paper we improve the
  known upper bounds on the numbers $F(2,2,p;p+1)$ and $F(3,p;p+1)$.

  \par
}

\renewcommand{\thefootnote}{}
\renewcommand{\thefootnote}{\arabic{footnote}}

\normalsize

\section{Introduction}

We consider only finite, non-oriented graphs without loops and multiple edges. We call a $p$-clique
of the graph G a set of $p$ vertices, each two of which are adjacent.  The largest positive integer
$p$, such that the graph G contains a $p$-clique is denoted by $cl(G)$ . We denote by $V(G)$ and
$E(G)$ the vertex set and the edge set of the graph $G $ respectively.  The symbol $K_n$ denotes
the complete graph on $n$ vertices.

Let $G_1$ and $G_2$ be two graphs without common vertices. We denote
by $G_1 + G_2$ the graph $G$ for which $V(G)=V(G_1) \cup V (G_2)$ and
$E(G)=E(G_1)\cup E(G_2) \cup E'$ , where $ E' = \{[x , y ] \mid x \in
V(G_1) , y \in V(G_2) \}$.

\begin{Definition}
  Let $ a_1 , \ldots , a_r $ be positive integers. We say that the
  $r$-coloring
$$
V(G)= V_1 \cup \ldots \cup V_r \text{ , } V_i \cap V_j = \emptyset
\text{ , } i \neq j,
$$
of the vertices of the graph $G$ is $(a_1 ,\ldots , a_r)$-free, if
$V_i$ does not contain an $a_i$-clique for each $i \in \{ 1 , \ldots ,
r \}$.  The symbol $G\rightarrow (a_1 , \ldots , a_r)$ means that
there is not an $(a_1 , \ldots , a_r)$-free coloring of the vertices
of $G.$
\end{Definition}

We consider for arbitrary natural numbers $a_1 , \ldots , a_r $ and
$q$
$$
H(a_1,\ldots , a_r ; q) = \{G : G\rightarrow (a_1,\ldots ,a_r) \text
{
  and } cl(G)<q \}.
$$
The vertex Folkman numbers are defined by the equality
$$ F(a_1, \ldots , a_r ;q) = \min\{ | V(G)| : G\in H(a_1, \ldots ,
   a_r; q)\}.
$$
It is clear that $G\rightarrow (a_1,\ldots ,a_r)$ implies $cl(G)\geq
\max\{a_1,\ldots , a_r\}$. Folkman [1] proved that there exists a
graph $G$ such that $G \rightarrow (a_1, \ldots , a_r) $ and $cl(G)
= \max\{a_1,\ldots , a_r\}$. Therefore
\begin{equation}
  \label{1} F(a_1 ,\ldots , a_r ;q) \text{ exists
    if and only if }
  q> \mbox{max} \{ a_1,\ldots , a_r \}.
\end{equation}
If $a_1 , \ldots ,a_r$ are positive integers, $r\geq2$ and $a_i=1$
then it is easy to see that
\begin{equation}
  G\rightarrow(a_1 , \ldots , a_r) \Leftrightarrow G \rightarrow(a_1 ,
  \ldots , a_{i-1} , a_{i+1} , \ldots ,a_r).
\end{equation}
It is also easy to see that for an arbitrary permutation $\varphi \in
S_r$ we have
$$
G\rightarrow(a_1, \ldots ,a_r) \Leftrightarrow
G\rightarrow(a_{\varphi(1)} , \ldots ,a_{\varphi(r)}).
$$
That is why
\begin{equation}
  F(a_1 , \ldots ,a_r;q)=F(a_{\varphi(1)} , \ldots , a_{\varphi(r)})
  ,\text{ for  each } \varphi \in S_r
\end{equation}
According to (2) and (3) it is enough to consider just such numbers $F(a_1, \ldots, a_r;q)$ for
which
\begin {equation}
  2 \leq a_1 \leq \ldots \leq a_r .
\end{equation}
For arbitrary positive integers $a_1, \ldots , a_r$ define:
\begin {gather}
  p=p(a_1 , \ldots ,a_r)= \max\{ a_1 ,\ldots , a_r \};\\
  m= 1+ \sum_{i=1}^{r} (a_i-1)
\end{gather}
It is easy to see that $ K_m\rightarrow(a_1 ,\ldots ,a_r)$ and
$K_{m-1}\nrightarrow(a_1 , \ldots ,a_r).$ Therefore
$$F(a_1 , \ldots , a_r;q)=m \text {, if  } q>m .       $$
In [4] it was proved that $F(a_1 ,\ldots , a_r;m)=m+p,$ where $m$ and $p$ are defined by the
equalities (5) and (6). About the numbers $F(a_1 , \ldots , a_r;m-1)$ we know that $F(a_1, \ldots ,
a_r;m-1) \geq m+p+2 \text { , } p \geq 2$ and
\begin{equation}
  F(a_1 , \ldots , a_r;m-1) \leq m+3p \text{ , [3]}.
\end{equation}

The exact values of all numbers $F(a_1 , \ldots , a_r;m-1)$ for
which $\max \{a_1 , \ldots , a_r \} \leq 4 $ are known. A detailed
exposition of these results was given in [8]. We must add the
equality $F(2,2,3;4)=14$ obtained in [2] to this exposition. We do
not know any exact values of $F(a_1 , \ldots , a_r;m-1)$ in the case
when $ \max \{a_1 , \ldots , a_r) \geq 5.$

According to (1), $F(a_1 , \ldots , a_r;m-1)$ exists exactly when $m \geq p+2.$ In this paper we
shall improve inequality (7) in the boundary case when $ m=p+2$ , $p \geq 5.$ From the equality
$m=p+2$ and (4) it easily follows that there are two such numbers only: $F(2,2,p;p+1)$ and
$F(3,p;p+1).$ It is clear that from $G\rightarrow (3,p)$ it follows $G\rightarrow(2,2,p).$
Therefore
\begin{equation}
  F(2,2,p;p+1) \leq F(3,p;p+1).
\end{equation}
The inequality (7) gives us that:
\begin{gather}
  F(3,p;p+1) \leq 4p+2;\\
  F(2,2,p:p+1) \leq 4p+2.
\end{gather}

Our goal is to improve the inequalities (9) and (10).  We shall need
the following

\begin{Lemma}
  Let $G_1$ and $G_2$ be two graphs such that
  \begin{equation}
    G_1 \rightarrow (a_1 , \ldots , a_{i-1}, a'_i , a_{i+1} , \ldots
    ,a_r)
  \end{equation}
  and
  \begin{equation}
    G_2 \rightarrow (a_1, \ldots ,a_{i-1} , a''_i , a_{i+1} , \ldots ,
    a_r).
  \end{equation}
  Then
  \begin{equation}
    G_1+G_2 \rightarrow (a_1 , \ldots , a_{i-1} , a'_i + a''_i ,  a_{i+1} , \ldots , a_r).
  \end{equation}
\end{Lemma}

\begin{proof}
 Assume that (13) is wrong and let
$$V_1 \cup \ldots \cup V_r \text { , } V_i \cap V_j = \emptyset
\text{ , } i \neq j$$

be a $(a_1 , \ldots ,a_{i-1} , a'_i + a''_i ,a_{i+1}, \ldots , a_r)-$free $r-$coloring of
$V(G_1+G_2)$.Let $V'_i = V_i \cap V(G_1) $ and $V''_i=V_i \cap V(G_2),$ for $i=1, \ldots , r.$ Then
$V'_1 \cup \ldots \cup V'_r$ is an $r-$coloring of $V(G_1),$ such that $V_j$ does not contain an
$a_j-$clique, $j \neq i.$ Thus from (11) it follows that $V'_i$ contains an $a'_i-$clique.
Analoguously from the $r-$colouring $V''_1 \cup \ldots \cup V''_r$ of $V(G_2)$ it follows that
$V''_i$ contains an $a''_i$-clique.Therefore $V_i=V'_i \cup V''_i$ contains a $(a'_i
+a''_i)-$clique,which contradicts the assumption that $V_1 \cup \ldots \cup V_r$ is a $(a_1, \ldots
, a_{i-1} , a'_i +a''_i , a_{i+1} , \ldots , a_r)-$ free $r-$coloring of $V(G_1+G_2).$ This
contradiction proves the Lemma.
\end{proof}

\section{Results}

The main result in this paper is the following

\begin{Theorem}
Let $a_1 \leq \ldots \leq a_r $ , $r \geq 2 $ be positive integers
and $ a_r = b_1 + \ldots + b_s,$ where $b_i$ are positive integers
,too and $b_i \geq a_{r-1}, \ $  $i=1 , \ldots , s.$ Then
\begin{equation}
  F(a_1 , \ldots , a_r; a_r+1) \leq \sum_{i=1}^{s}
  F(a_1 , \ldots , a_{r-1} , b_i;b_i +1).
\end{equation}
\end{Theorem}

\begin{proof}
 We shall prove the Theorem by induction on $s.$ As
the inductive step is trivial we shall just prove the inductive base
$s=2.$ Let $G_1$ and $G_2$ be two graphs such that $cl(G_1)=b_1$ and
$ cl(G_2)=b_2,$ $a_r =  b_1 +b_2 , \  b_1 \geq a_{r-1}, \ b_2 \geq
a_{r-1} $ and
\begin{gather*}
G_1 \rightarrow (a_1 , \ldots , a_{r-1} , b_1) \text{ , }
 |V(G_1)|= F(a_1 , \ldots , a_{r-1},b_1;b_1 + 1)\\
G_2 \rightarrow (a_1 , \ldots , a_{r-1}, b_2) \text{ , } |V(G_2)|
= F(a_1 , \ldots , a_{r-1} , b_2 ;b_2 + 1).
\end{gather*}
According to the Lemma, $G_1+G_2 \rightarrow (a_1, \ldots , a_{r-1} , a_r).$ As $cl(G_1 + G_2)
=cl(G_1) + cl(G_2) =b_1+ b_2 = a_r, $ we have
$$F( a_1 , \ldots , a_r;a_r +1) \leq |V(G_1 + G_2)| = |V(G_1)| +
|V(G_2)|.$$

From this inequality (14) trivially follows when $s=2$ and hence , for arbitrary $s,$ as explained
above. The Theorem is proved.
\end{proof}

We shall derive some corollaries from the Theorem. Let $p \geq 4$
and $p=4k+l, \ $    $0 \leq l \leq 3.$ Then from (14) it easily
follows that
\begin{gather}
  F(3,p;p+1) \leq (k-1) F(3,4,;5) + F(3, 4+l;5+l) \\
  F(2,2,p; p+1) \leq (k-1) F(2,2,4;5) + F(2,2,4 +l ;5+l).
\end{gather}
From (15), (9) $(p=5 , 6 , 7)$ and the equality $F(3,4;5) = 13,$ [6]
we obtain

\begin{corollary}
Let $p \geq 4$. Then:
\begin{gather*}
  F(3,p;p+1) \leq \frac{13p}{4}  \  \text{ for }  \ p \equiv 0 \mod 4; \\
  F(3,p;p+1) \leq \frac {13p+23}{4} \  \text { for } \ p \equiv 1 \mod 4;\\
  F(3,p;p+1) \leq \frac {13p+26}{4} \  \text { for } \ p \equiv 2 \mod 4;\\
  F(3,p;p+1) \leq \frac {13p+29}{4} \  \text { for } \ p \equiv 3 \mod 4.
\end{gather*}
\end{corollary}

From (16), the equality $F(2, 2 ,4;5)=13,$ [7], the inequality (10)
$(p=5)$ and the inequalities $F(2, 2 , 6;7 ) \leq 22,$ [9];
$F(2,2,7;8) \leq 28 ,$ [9] we obtain

\begin{corollary}
Let $p \geq 4$. Then
\begin{gather*}
  F(2,2,p;p+1) \leq \frac {13p}{4} \  \text { for } \ p \equiv 0 \mod 4;\\
  F(2,2,p;p+1) \leq \frac {13p+23}{4} \  \text { for } \  p \equiv 1 \mod 4;\\
  F(2,2,p;p+1) \leq \frac {13p+10}{4}  \ \text { for } \ p \equiv 2 \mod 4;\\
  F(2,2,p;p+1) \leq \frac {13p+21}{4} \  \text { for } \ p \equiv 3 \mod 4.
\end{gather*}
\end{corollary}

We conjecture that the following inequalities hold:
\begin{gather}
  F(3,p;p+1) \leq  \frac {13p}{4} \  \text {   for   } p \geq 4 ;\\
  F(2,2,p;p+1) \leq  \frac {13p}{4}  \ \text {  for  } p \geq 4 .
\end{gather}

From the Theorem it follows that
\begin {equation}
  F(3,p;p+1) \leq F(3,p-4 ; p-3) + F(3,4;5) ,p \geq 8 ;
\end{equation}

\begin {equation}
  F(2,2,p;p+1) \leq F(2,2,p-4 ;p-3) + F(2,2,4;5) , p \geq 8.
\end{equation}

From  $F(3,4;5) = 13$ (see [6]) and (19) we obtain

\begin{corollary}
If the inequality (17) holds for $p=5$, $6$ and $7$, then (17) is true for
every $p \geq 4$.
\end{corollary}

From $F(2 , 2 , 4 ;5) = 13,$ (see [7]) and from (20) it follows

\begin{corollary}
If the inequality (18) holds for $p=5$, $6$
  and $7$ then (18) is true for every $p \geq 4$.
\end{corollary}

At the end in regard with (8) we shall pose the following

\begin{Problem}
Is there a positive integer $p,$ for which $F(2,2,p;p+1) \neq
F(3,p;p+1)$?
\end{Problem}

\bigskip \medskip

\parbox[t]{6cm}{
  \normalsize
  N. Kolev\\
  \small
  Department of Algebra,\\
  Faculty of Mathematics and Informatics,\\
  ``St. Kl. Ohridski'' University of Sofia,\\
  5 J. Bourchier blvd, 1164 Sofia,\\
  BULGARIA\\
  \texttt{nickyxy@mail.bg}
}
\hfill
\parbox[t]{6cm}{
  \normalsize
  N. Nenov\\
  \small
  Department of Algebra,\\
  Faculty of Mathematics and Informatics,\\
  ``St. Kl. Ohridski'' University of Sofia,\\
  5 J. Bourchier blvd, 1164 Sofia,\\
  BULGARIA \\
  \texttt{nenov@fmi.uni-sofia.bg}
 }


\end{document}